# Sequential Hybrid Finite Element and Material Point Method to Simulate Slope Failures

Brent Sordo[1], Ellen Rathje, Ph.D.[2], Krishna Kumar, Ph.D.[3]

[1] Ph.D. Candidate, Department of Civil, Architectural, and Environmental Engineering, The University of Texas at Austin, 301 E Dean Keeton St, Austin, TX 78712; E-mail: bsordo@utexas.edu

[2] Professor, Department of Civil, Architectural, and Environmental Engineering, The University of Texas at Austin, 301 E Dean Keeton St, Austin, TX 78712; E-mail: e.rathje@mail.utexas.edu

[3] Assistant Professor, Department of Civil, Architectural, and Environmental Engineering, The University of Texas at Austin, 301 E Dean Keeton St, Austin, TX 78712; E-mail: krishnak@utexas.edu

**ABSTRACT**

Numerical modeling of slope failures seeks to predict two key phenomena: the initiation of failure and the post-failure runout. Currently, most modeling methods for slope failure analysis excel at one of these two but are deficient in the other. For example, the Finite Element Method (FEM) models the initiation of instability well but quickly loses accuracy when modeling large deformations because of mesh distortion, restricting its ability to predict runout. Conversely, the Material Point Method (MPM) utilizes material points which move freely across a background grid, allowing for indefinite deformations without computational issues. However, MPM is restricted in its ability to model slope failure initiation due to limitations of the available boundary conditions and reduced accuracy of its stress distributions. The sequential hybridization of these two methods, initiating a model in FEM and then transferring to MPM, presents an opportunity to accurately capture both initiation and runout by a single model. The exact time for this transfer is not self-apparent, but it must be conducted after the initiation mechanism and before excessive mesh distortion. By simulating two granular column failures and two slope failures, we demonstrate the effectiveness of this hybrid FEM-MPM method and identify the appropriate time to transfer.

# INTRODUCTION

Landslides represent major geohazards arising from both natural and anthropogenic factors. They involve rapid down-slope movements of soil and rock, which can cause extensive damage to infrastructure and pose severe threats to human safety. In response to these challenges, geotechnical engineers face the critical task of forecasting these events to support risk assessment, to develop strategies for avoidance, and to formulate mitigation measures. Numerical modeling has emerged as an essential tool in assessing slope vulnerabilities and informing these vital decisions. Effective slope stability analysis requires a comprehensive approach, encompassing two fundamental aspects: (1) determining the specific conditions that trigger a slope failure and (2) estimating the magnitude and trajectory of landslide movements in the event of a failure.

Traditional slope stability analysis methods, such as the Limit Equilibrium Method (LEM), only satisfy stress equilibrium with a rigid, perfectly-plastic failure criterion and ignore strain compatibility (Huang, 2014). In contrast, the Finite Element Method (FEM) model of a slope failure considers stress equilibrium, strain compatibility, and a comprehensive constitutive law (Zienkiewicz and Taylor, 1967). Thus, FEM can accurately identify conditions which cause the failure of a slope and the initial development of the failure surface (Naylor, 1982), even for complex geometries (Zheng et al., 2005; Zheng et al., 2008). Nevertheless, FEM has limitations in modeling large deformations. While adept at analyzing small strain problems, such as the initiation of slope failures, FEM loses accuracy when predicting extensive landslide runout due to mesh distortion and entanglement under large deformations. While re-meshing (Zhang et al., 2016) and displacement field refinement (Zeng and Liu, 2016) mitigate mesh distortion issues, they require substantial computational power. This makes them impractical for modeling large landslide runouts (Liang and Zhao, 2018).

Consequently, the inherent limitations of FEM restrict its ability to reliably forecast landslide runouts and their consequent downslope ramifications (Soga et al., 2015). The Material Point Method (MPM; Sulsky et al., 1994; Bardenhagen et al., 2000) is an alternative numerical approach suitable for modeling large deformation problems. MPM is a hybrid Eulerian-Lagrangian approach. The domain is modeled as individual material points that carry historic state parameters and can move freely through a background mesh and deform according to Newtonian laws without distorting the background mesh. Thus, MPM can account for large displacements, such as landslide runouts, without suffering from mesh distortions (Cuomo et al., 2021). However, MPM integrates its stresses at the material points. As material points traverse the mesh, the integration accuracy diminishes because the integration weights stay fixed, despite the shift in the integration location (i.e., the position of the material points). This issue contrasts with the Gauss integration method used in the FEM, where the integration is always done at specific Gauss locations and appropriate Gauss weights. Furthermore, as material points crosses from one grid cell to another, it produce spurious stress oscillations leading to significant errors, termed as cell crossing issue (Bardenhagen et al., 2000). This material point integration results in severe stress checkerboarding (Wang et al., 2021). Due to its dual representation system, MPM also faces challenges modeling certain kinematic boundary conditions like absorbing boundaries and free-field columns for earthquake analysis; particles carrying information about the material state can move through the background mesh, but boundary conditions are applied at the mesh nodes. While dashpots and damping layers have been utilized in MPM (Shan et al., 2021), the potential for material points to move away from the boundary cells restricts its options for boundary conditions. These stress inaccuracies and boundary condition limitations restrict the ability of MPM to model complex initiation scenarios such as earthquake-induced liquefaction. Although MPM can model large deformations, its shortcomings associated with cell crossing, stress checkerboarding, and free-field boundary conditions limit its accuracy in predicting failure initiation.

A comprehensive slope failure analysis requires accurate failure initiation and subsequent runout modeling. Stresses and boundary conditions determine critical features of the failure initiation, such as the geometry of the failure surface, the volume of the sliding mass, and the pore pressure response (Federico et al., 2015). Meanwhile, the initiation energy and the topographical features largely govern the runout response.

To overcome the shortcomings of each method, we propose a sequential hybrid FEM-MPM method to simulate the entire landslide process from initiation using FEM to runout using MPM. The hybrid approach combines the advantages of both numerical techniques to comprehensively evaluate the entire process of complex landslide failures, including those induced by earthquake-shaking and liquefaction, which neither method would be capable of independently. To our knowledge, this is the first time FEM and MPM have been combined sequentially to fully utilize the strengths of both methods over the same soil mass, although they have been combined in other ways.

FEM and MPM have been utilized simultaneously within the same model, with each simulating a different section of the domain (Lian et al., 2012; Pan et al., 2021). These studies focused on the interaction of soft and stiff materials in which MPM models the soft material that experiences large deformation, and the stiff material is modeled by FEM characterized by its small deformations. FEM has also been combined with the smoothed particle hydrodynamic method (SPH; Groenenboom et al., 2019; Long et al., 2020) to simulate large deformation fluid-structure interaction problems. In these studies, FEM is used to simulate the solid domain, and SPH models the fluid domain.

While useful in their applications (i.e., the interaction of solids and fluids or soft and stiff materials), these hybrid approaches are focused on spatial hybridization rather than the temporal or sequential hybridization required to simulate earthquake-induced landslides, one of the main motivations of our hybrid method. Predicting the triggering of liquefaction in a slope requires accurate stresses, an advanced constitutive model, and dynamic boundary conditions, all of which necessitate the use of the FEM. After liquefaction triggering, the soil experiences significant strength reduction and thus transitions to a highly deformable material, requiring MPM simulations. Thus, to capture the response of earthquake-induced failures, we propose to employ FEM and MPM sequentially over time across the entire domain to capture the failure initiation when the soil experiences small-strain behavior and the post-liquefaction runout when the soil undergoes large deformation. No current methods exist which hybridize FEM and MPM sequentially over the same domain to capture the consequences of this fundamental transformation of material behavior.

Researchers have actively developed hierarchical multiscale hybridizations, building upon the need for hybrid methods to model transitions in material behavior from solid-like to fluid-like responses. Specifically, Liang and Zhao (2018) combined the MPM with the Discrete Element Method (DEM), enabling the constitutive stress response within the MPM algorithm to be determined by micro-scale DEM models of soil grains at material point contacts. Andrade and Mital (2019) integrated DEM with the FEM and the Finite Difference Method (FDM) in a similar manner. These hybrid methods have significantly enhanced the accuracy of the stress response by capturing nonlinearity and state dependence in the constitutive models. However, they entail using continuum and discrete frameworks within a single simulation, leading to computational inefficiency (Liang and Zhao, 2018). Moreover, these latter methods have not addressed the limitations of FEM in handling large deformation runouts. As a result, despite these advancements, multiscale hybridizations remain inefficient at the prediction of earthquake-induced liquefaction and subsequent complex slope failure.

Nevertheless, there have been some sequential hybridizations of FEM with SPH that have yielded success. Kitano et al. (2017) created models of embankment slope failures in which finite elements were converted into SPH particles when the elements exceeded a shear threshold. However, these models were primarily concerned with determining the impact of the slope failure on embedded pipelines, not the runout of debris. Additional research concerning sequential hybridizations for the application of earthquake-triggered landslide runouts is therefore necessary.

Talbot et al. (2024) developed a sequential hybrid FDM-MPM model to capture the large deformation failure of the San Fernando Dam due to earthquake-induced liquefaction. Their hybrid FDM-MPM method creates MPM particles at the Gauss point locations of the FD elements which directly inherit the stresses of that respective FD Gauss point. They also designate certain MPM particles to have softened material properties if they are located in the strain bands in the FDM or are liquefied. However, they do not quantitatively transfer strains, velocities, and other state variables from FDM to the material points. Because of this, their method is limited to cases in which all kinematics develop in the MPM phase. The San Fernando Dam failure uniquely fits this criterion as the dam collapse did not begin until after shaking had finished. Thus, their hybrid model with MPM points based only on the FDM stresses predicts collapse similar to the actual failure. Their model demonstrates the viability of sequentially hybridizing Lagrangian and Eulerian-Lagrangian methods, but their particular method is limited to transfers where the kinematics only exist in the MPM phase.

Our proposed hybrid FEM-MPM method involves transferring the entire problem domain, including all its properties and state variables, from FEM to MPM as the failure evolves from the initiation to the runout stage. Thus, our method includes the transfer of kinematics from FEM to MPM, making it applicable in scenarios where the failure initiates during earthquake shaking. Furthermore, it interpolates state variables throughout the model domain, allowing for the creation of any number of particles in the MPM phase and giving the user more options regarding model resolution.

We investigate the influence of the transfer time ($t_T$) on the failure kinematics and runout of simple, gravity-driven failures: granular column collapse experiments and slope failures. The granular column collapse (Lajeunesse et al., 2005; Soundararajan, 2015) is a simple and well-documented experiment of granular flow dynamics that captures the dynamics of large landslides (Lajeunesse et al., 2005). We assess how transfer time affects model performance by evaluating multiple sequential hybrid FEM-MPM models at different transfer times. We also establish criteria linking the accuracy of the runout profile to the characteristics of the FEM mesh at the transfer time. While these simple, gravity-driven failure simulations do not fully utilize the potential of the hybrid FEM-MPM method as they do not involve complex failure initiation mechanisms, such as earthquake-induced liquefaction. However, they are valuable for demonstrating the effectiveness of the hybrid method, assessing how the transfer time ($t_T$) influences the movement and final spread of the failure, and guiding the selection of suitable transfer times for more intricate failure scenarios.

**HYBRID FEM-MPM APPROACH**

The proposed sequential hybrid FEM-MPM approach consists of three phases (Figure 1): (a) an initial FEM phase to capture failure initiation, (b) a transfer phase that moves the data from FEM to MPM, and (c) a final runout phase in MPM. The hybrid approach starts by initializing the model in FEM to a geostatic stress state, employing linear elastic constitutive models for all materials. During the failure initiation phase, FEM switches from linear elastic materials for geostatic stresses to appropriate plastic constitutive models. This initial FEM phase progresses until a user-specified transition time ($t_T$), at which time we perform the transfer. The transfer phase consists of moving the geometry and the properties from FEM to MPM. We transfer these properties from FEM to material points by generating a user-specified number of material points for each FEM element (Figure 2), via a Python script. These material points carry stress, kinematic information, and other relevant state variables from the FEM phase. We develop a generic transfer algorithm that can handle any type of FEM element and any number of material points per element. This algorithm assigns material points to Gauss locations within each element, according to the user-defined number of Gauss locations. However, the specific locations and quantity of Gauss points utilized during the Finite Element Method (FEM) phase might differ from those selected for initializing the material points.

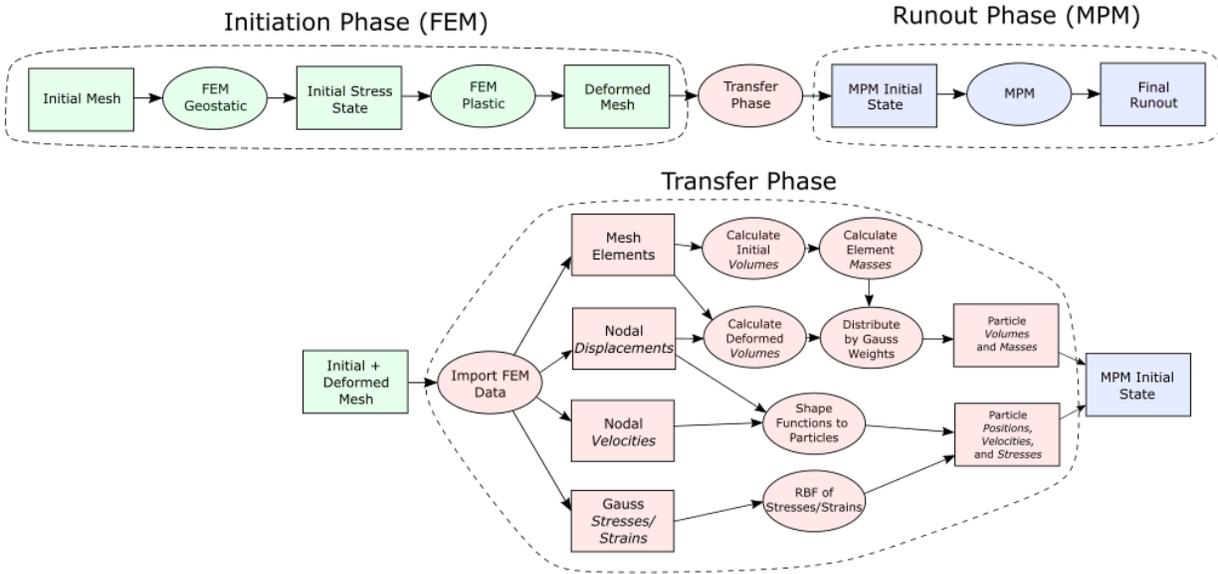

*Figure 1: Workflow of sequential hybrid FEM-MPM method*

In FEM, stresses and strains are stored at the Gauss locations, while kinematic information, such as displacements and velocities, is recorded at the nodes. We use isoparametric shape functions to interpolate state variables from the FEM nodes to the material points. To interpolate stresses from FEM Gauss points to the material points, we employ a radial basis function (Travers, 2007). We assign volume and mass of the material points based on the geometry of the elements in the FEM phase. The MPM phase utilizes an independent background mesh that is different from the FEM mesh. We initialize the MPM phase with material points carrying information about positions, velocities, pore pressures, stresses, strains, volume, and mass determined by the transfer algorithm and then run the MPM phase until the runout process is complete. All transfer properties are stored in the material points, and the background mesh does not carry any information at the transfer time.

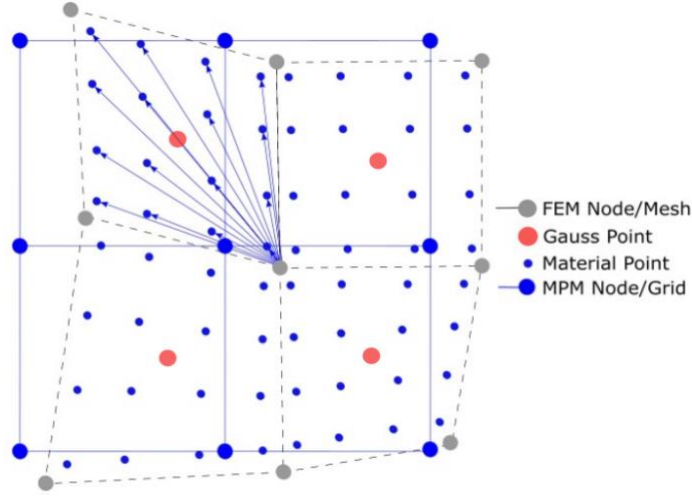

*Figure 2: Schematic of FEM to MPM transfer procedure illustrating the interpolation of state variables from nodes to material points. This example utilizes one Gauss point per element and sixteen material points per element.*

*Transfer Formulation*

The FEM nodes store information about geometry, displacement, velocities, and pore pressures for each element. We employ $C^0$ isoparametric shape functions for each FEM element to transfer these nodal quantities to the material point locations, as shown in Figure 2, according to:

$$u^P(\xi_i, \eta_j) = \sum_{k=1}^{N_n} N_k(\xi_i, \eta_j) * u_k \quad (1)$$

where $N_n$ denotes the total number of nodes, $u^P$ is the state variable assigned to the particle, $\xi_i$ and $\eta_j$ are the isoparametric coordinates of the Gauss location of the particle, $N_k$ is the shape function for node $k$, and $u_k$ is the value of the state variable at node $k$.

In FEM analysis, the stresses and strains are stored at the Gauss points. Because the FEM Gauss points are not co-located with the material points, it is possible that equilibrium will not be maintained in the MPM model during interpolation. To minimize stress disequilibrium, we interpolate stresses and strains from the FEM Gauss points to the material points using a multi-quadratic radial basis function (RBF; Travers, 2007). The fundamental concept behind RBF is to construct a smooth function that exactly passes through a given set of points. This method is especially useful in scattered data, where the data points are irregularly distributed in space, such as a deformed FEM mesh. The multi-quadratic RBF for a set of points is defined as:

$$\Phi(r) = \sqrt{(r/c)^2 + 1} \tag{2}$$

where $r$ is the Euclidean distance between any two points in space, $r = |x - x_i|$, $c$ is a positive constant known as the shape parameter, $x$ is the point where the function is evaluated (i.e., material point), and $x_i$ is a center of the RBF (i.e., Gauss points in FE). We assign $c$ the default value recommended by Travers (2007), which is the average distance between the FEM Gauss points.

Given a set of $N_g$ FEM Gauss points $\{(x_i, u_i)\}_{i=1}^{N_g}$ where $x_i$ is the Gauss locations in space and $u_i$ are the stress and strain values associated with these Gauss points, the goal is to find an interpolating function, $s(x)$, that exactly fits these points. This function is constructed as a linear combination of RBFs centered at each data point:

$$s(x) = \sum_{i=1}^{N_g} \lambda_i \Phi(|x - x_i|) \tag{3}$$

where $\lambda_i$ are coefficients to be determined.

The coefficients $\lambda_i$ are found by solving the linear system formed by enforcing the interpolation conditions $s(x_i) = y_i \; for \; all \; i$. This results in a system of linear equations:

$$\begin{bmatrix} \phi(\|x_1 - x_1\|) & \cdots & \phi(\|x_1 - x_N\|) \\ \vdots & \ddots & \vdots \\ \phi(\|x_N - x_1\|) & \cdots & \phi(\|x_N - x_N\|) \end{bmatrix} \begin{bmatrix} \lambda_1 \\ \vdots \\ \lambda_N \end{bmatrix} = \begin{bmatrix} y_1 \\ \vdots \\ y_N \end{bmatrix} \tag{4}$$

RBF can also extrapolate the material points near the model boundaries via the same function. Once the interpolating function $s(x)$ is constructed and the coefficients $\lambda_i$ are determined, this function can be evaluated at any point $x$ in space, which allows for extrapolation smoothly and continuously based on the behavior of the data within its range.

To evaluate the stress equilibrium during transfer using RBF, we consider transferring the vertical stress of a 4 m x 4 m geo-static linear elastic column of dry soil. The FEM analysis consists of 256 square elements with an element size of 0.25 m, each with a single, center Gauss point, and each FEM element consists of 16 material points. Figure 3a shows the vertical stress distribution at the material points after RBF interpolation from the FEM simulation, while Figure 3b shows the error in the RBF stresses at the material points relative to the theoretical vertical stress values calculated for geostatic conditions. The multi-quadratic RBF smoothly interpolates stresses within the inner domain with essentially no error, but close to the boundaries the RBF stresses have a small error due to extrapolation. These small errors are considered acceptable and we do not believe they will significantly influence the final results. This issue is explored further later in the paper by comparing runout results between a pure-MPM model and a hybrid FEM-MPM model.

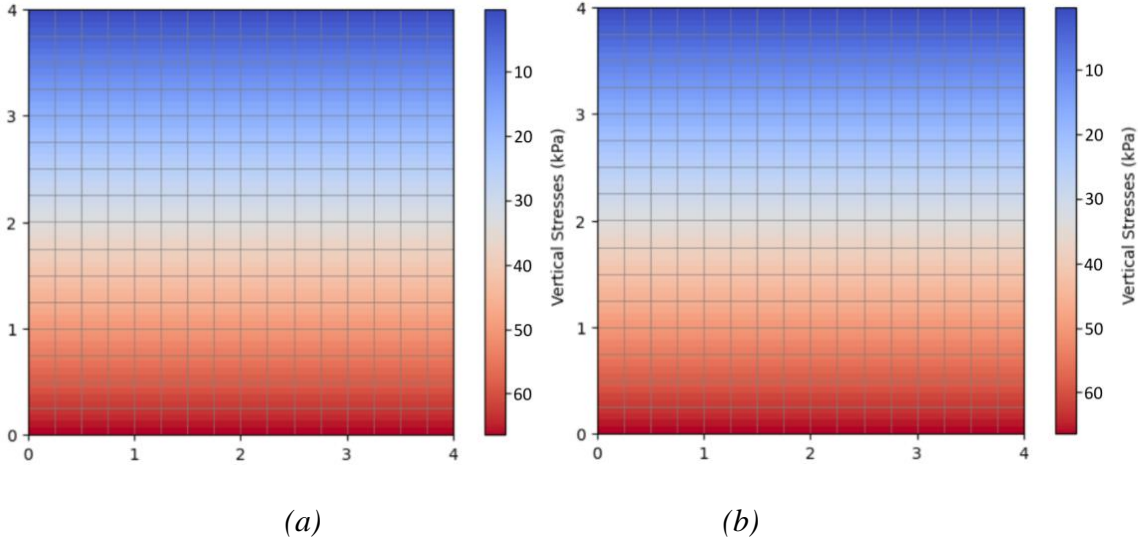

*(a)*                      *(b)*

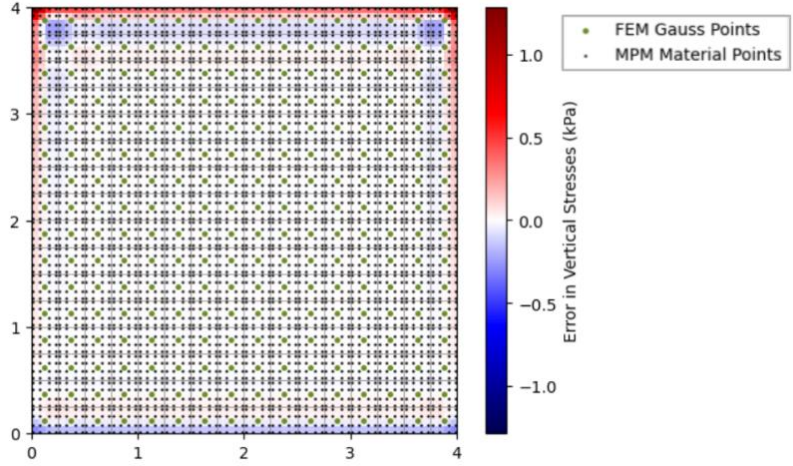

*(c)*

*Figure 3: (a) Theoretical stress distribution in a 4 m x 4 m linear elastic column after transfer from FEM. (b) Stresses resulting from RBF interpolation. (c) Error in vertical stresses after the transfer using RBF interpolation with 0.25 m x 0.25 m elements.*

It is important to accurately transfer the FEM element volume and mass to the material points, conserving volume and mass during the transfer. In this study, we transfer the volume of the deformed elements at the time of the transfer and the mass of the undeformed elements. The volume of a deformed FEM element at the transfer time $t_T$ ($V^e_t$) is transferred to ($N_{px}$ by $N_{py}$) material points within that element using the Gauss-Legendre weights at the material point locations.

$$V^e_t = \sum_{i=1}^{N_{px}} \sum_{j=1}^{N_{py}} V^P_{i,j} \qquad (5)$$

and

$$V^P_{i,j} = V^e_t * w(\xi_i)/2 * w(\eta_j)/2 \qquad (6)$$

where $V^P_{i,j}$ is the volume assigned to the particle located at isoparametric location $\xi_i$ *and* $\eta_j$, i and j take on values [1, $N_{px}$] and [1, $N_{py}$], respectively, and $w(\xi_i)$ and $w(\eta_j)$ are the Gauss weights at locations $\xi_i$ *and* $\eta_j$, each normalized by 2 to ensure the sum of the Gauss-Legendre weights equal to 1.0. As a result, the sum of $V_{i,j}^P$ for equals the volume of the element ($V^e{}_t$). To conserve the mass of the original FEM phase, we assign mass to each material point using the initial volume of each FEM element ($V_0^e$) and its density:

$$M^e = \sum_{i=1}^{N_{px}} \sum_{j=1}^{N_{py}} M_{i,j}^p \tag{7}$$

and

$$M_{i,j}^P = V_0^e * \rho * w(\xi_i)/2 * w(\eta_j)/2 \tag{8}$$

where $M^e$ is the mass of the element, $M_{i,j}^P$ is the mass assigned to the material point at location $\xi_i$ and $\eta_j$, $V_0^e$ is the initial volume of the element, $\rho$ is the mass density of the material, and $w(\xi_i)$ and $w(\eta_j)$ are the Gauss weights at location $\xi_i$ and $\eta_j$.

*MPM Formulation*

MPM segments the domain, $\mathcal{B}$, into $N_p$ material points, each with a volume of $V^p$:

$$\mathcal{B} = \sum_{p=1}^{N_p} V_p \tag{9}$$

Applying the balance of linear momentum, the weak form in the updated Lagrangian framework for a body, $\mathcal{B}$, over the material domain, $\Omega$, is written in the discretized form as:

$$\sum_{i,j=1}^{n_n} \int_\Omega \rho N_i(\mathbf{x}) N_j(\mathbf{x}) \mathbf{a}_j(t) \, d\Omega = -\sum_{i=1}^{n_n} \int_\Omega \nabla N_i(\mathbf{x}) : \boldsymbol{\sigma}(\mathbf{x},t) \, d\Omega$$

$$+ \sum_{i=1}^{n_n} \int_{\partial \Omega_\Gamma} N_i(\mathbf{x}) \mathbf{t}(\mathbf{x},t) \, dS$$

$$+ \sum_{i=1}^{n_n} \int_\Omega \rho N_i(\mathbf{x}) \mathbf{b}(\mathbf{x},t) \, d\Omega \tag{10}$$

where $M^p$ denotes the mass of a particle $p$. The subscript $p$ indicates that a term describes a variable associated with the material point $p$.

MPM considers the locations of material particles as the integration points. This converts the integrals in the weak form into sums over the material particles. As a result, the semi-discrete momentum equation can be written as:

$$\sum_{j=1}^{N_n} M_{ij} a_j = F_i^{int} + F_i^{ext} \qquad (11)$$

where $N_n$ is the number of nodes, $M_{ij}$ is the consistent mass matrix, $F_i^{int}$ is the internal force vector and $F_i^{ext}$ is the external force vector.

In the standard MPM formulation, the consistent mass matrix is generally replaced by the lumped masses at nodes, $m_i$. The internal force vector computed at the grid nodes is written as

$$F_i^{int} = -\sum_{p=1}^{N_p} \nabla N_i(x^p) : \boldsymbol{\sigma}(x^p, t) V^p \qquad (12)$$

and the external force vector is:

$$F_i^{ext} = \int_{\partial \Omega_\Gamma} N_i(x) \boldsymbol{t}(x, t) dS + \sum_{p=1}^{N_p} N_i(x^p) \boldsymbol{b}(x^p, t) M^p \qquad (13)$$

MPM requires mapping quantities between nodes and material points. We transfer quantities between material points and mesh nodes using shape functions. For example, the displacement of a material point can be determined by interpolating the displacements at the nodes:

$$u_p = \sum_{i=1}^{N_n} N_i^p u_i \qquad (14)$$

where $N_n$ denotes the total number of nodes, and $N_i^p$, which is equivalent to $N_I(\xi_p)$, signifies the value of the shape function for node $I$ at material point $p$, with $\xi_p$ being the local coordinates of the material point. Similarly, other material point quantities, such as velocity and acceleration, are interpolated from their respective values at the grid nodes. For more information on MPM, see Soga et al., 2016.

# EVALUATING THE SEQUENTIAL HYBRID FEM-MPM APPROACH USING GRANULAR COLUMN COLLAPSE EXPERIMENT

## *Granular Column Collapse*

Granular column collapse is a simple experiment that captures the dynamics of granular flow behavior, such as landslides (Lajunesse et al., 2005). The experiment involves the sudden collapse of a granular soil column of initial height ($H_0$) and length ($L_0$) onto a horizontal surface. The collapsing soil flows freely until the movement ceases due to friction between the soil and the horizontal surface, generating a final runout geometry. The runout is defined as the distance the toe moves from its original location, and normalized runout ($R_N$; Kumar and Soga, 2019) is defined as:

$$R_N = \frac{R}{L_0} = \frac{L - L_0}{L_0} \tag{15}$$

where L and $L_0$ represent the current and initial column widths, respectively.

The collapse occurs in three phases: initial acceleration, runout, and cessation (Figure 4). Lajeunesse (2005) introduced the concept of a critical time ($\tau_c$), which represents the duration of the initial acceleration phase as:

$$\text{critical time} = \tau_c = \sqrt{\frac{H_0}{g}} \tag{16}$$

where *g* is the acceleration of gravity. During the initial acceleration phase (Figure 4), a failure plane develops from the toe of the column, causing a collapse along the flank. The material above the failure plane experiences a sliding failure for short columns and a free fall under gravity for tall columns. This initial acceleration phase extends until the collapse reaches its peak kinetic energy, termed "full mobilization," at *t / $\tau_c$* ~ 1, after which the runout phase begins (Figure 4). The runout phase is characterized by the transition of vertical acceleration to horizontal acceleration as the materials interact with the horizontal surface. At $t/\tau_c \sim 3$, the cessation phase begins (Figure 4) as the runout nears completion. The movement decelerates due to frictional resistance at the base, and the flow reaches its final geometry by $t/\tau_c \sim 5$.

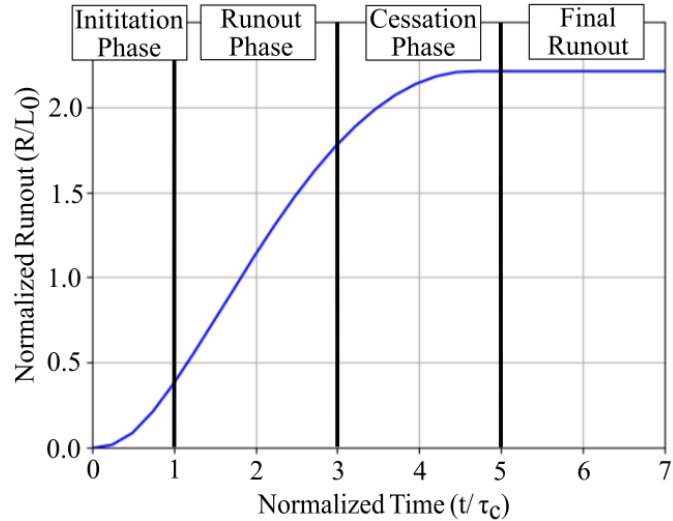

*Figure 4: Schematic of normalized runout over time for a pure MPM column collapse demonstrating the three collapse phases.*

We demonstrate the performance of the hybrid FEM-MPM method in predicting the runout of a short granular column with an aspect ratio, $a$, of 0.8 ($a = H_0/L_0 = 40\ m/50\ m$). For this example, we transfer the geometry and material state of the FEM phase to MPM at the critical time ($t = \tau_c = \sqrt{40\text{m}/9.81\text{m}/\text{s}^2} \sim 2.0$ s) of the column and compare these results to pure MPM and pure FEM analysis.

For the FEM phase, we use the open-source FEM code OpenSees (McKenna, 1997) and discretize the granular column with 1 m x 1 m elements. We use SSPquadUP elements (McGann et al., 2012), which are reduced-order, linear, quadrilateral elements with four nodes and a single central Gauss point. We simulate the geostatic stage of the column using a linear elastic model with properties shown in Table 1. The model initially includes roller boundaries on the sides and a frictional boundary on the bottom with a coefficient of friction of 0.466, matching the experiments of Lajeunesse (2005). After the initial geo-static stress state, we initiate failure by removing the roller boundary on the right side and switching to a plastic constitutive model with the properties defined in Table 1. We specify the material parameters of the nonlinear Pressure Independent Multi-Yield (PIMY; Yang et al., 2005) constitutive model to represent the Mohr-Coulomb failure criteria, as described in Hwang and Rathje (2023). The elastic properties remain unchanged between the geostatic and plastic failure phases. We use the strength properties reported from the experimental results of Lajeunesse (2005) and include 1 kPa of cohesion to ensure numerical stability. For the MPM runout phase, we employ the CB-Geo MPM code (Kumar et al., 2019) with the Generalized Interpolation Material Point (GIMP) method (Bardenhagen et al., 2004). The MPM model uses a square background grid of 1 m x 1 m cells. For each FEM element, we generate 16 material points. The column is modeled using the Mohr-Coulomb failure criterion with material properties matching its FEM counterpart.

*Table 1: Properties of granular column material*

| Mass density | 1700 | kg/m3 |
|---|---|---|
| **Elastic parameters** | | |
| Young's modulus | 23.8 | MPa |
| Poisson's ratio | 0.23 | |
| **Strength parameters** | | |
| Friction angle | 22.2 | degrees |
| Cohesion | 1 | kPa |

Figure 5 shows the geometries and distributions of stresses at different stages of the sequential hybrid FEM to MPM analysis. Figures 5a and 5b show the vertical stress distribution immediately before (FEM) and after (MPM) the transfer time ($t_T = \tau_c = 2.0$ s). While there are minor differences in the stresses after the transfer, the transfer algorithm generally preserves stress equilibrium. Figure 5c shows the vertical stress distribution at the MPM material points at $t = 2.5$ s, 0.5 s into the MPM phase. At this point, the MPM stresses exhibit severe checkerboarding while the kinematics remain unaffected. Finally, Figure 5d shows the MPM results for the final runout ($t = 20.0$ s). The MPM phase delivers a smooth, realistic runout profile, but stress checkerboarding persists even after movement has ceased.

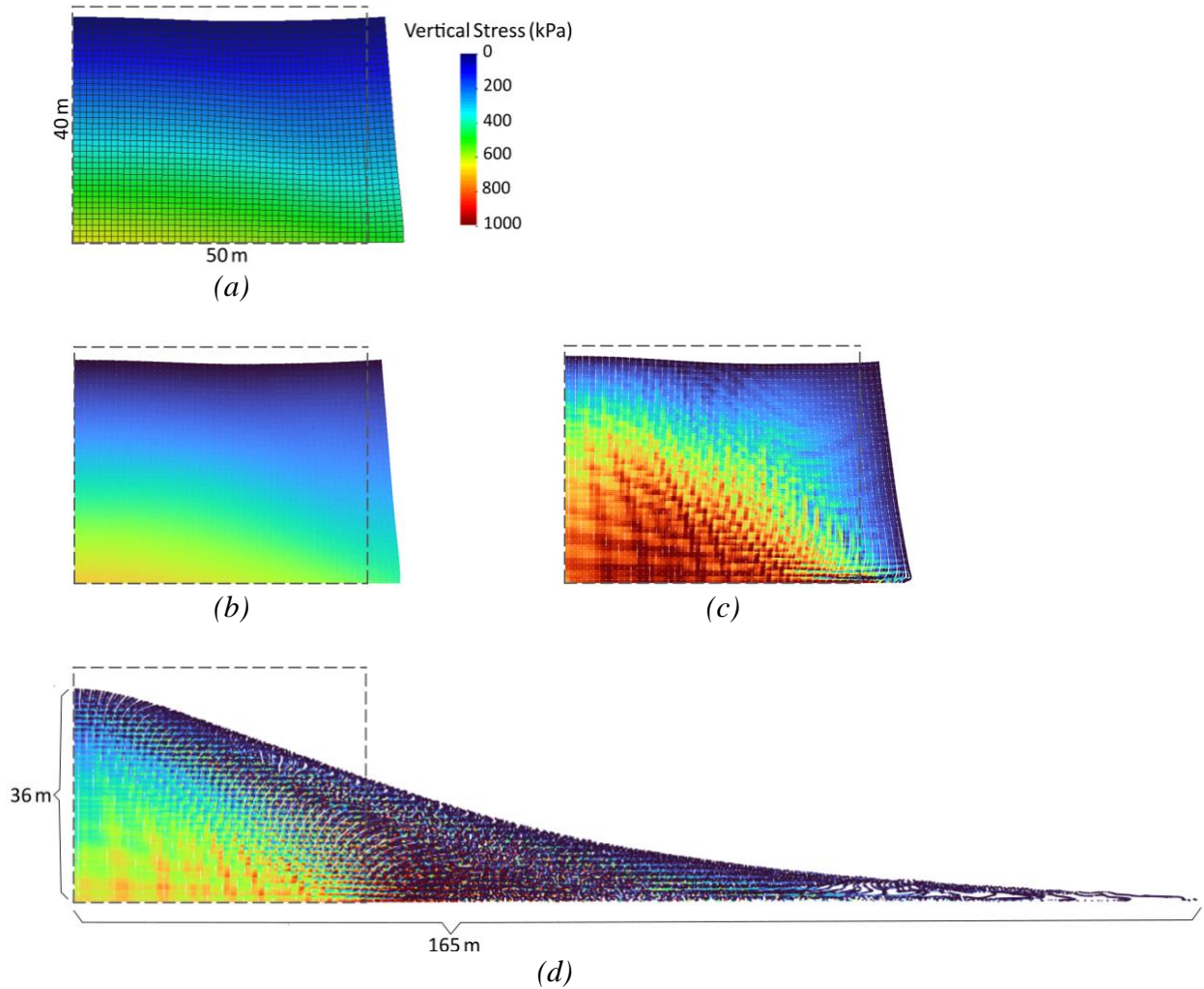

*Figure 5: Results from hybrid FEM-MPM model with transfer at t = 2.0 s. Stress contours from a) FEM at transfer, b) MPM material points immediately after transfer, c) MPM material points 0.5 s after transfer (t = 2.5 s), and c) MPM material points at the end of runout (t = 20.0 s).*

Figure 6 shows the final runout profile of the granular column collapse from our hybrid FEM-MPM column model transferred at $t_T = 2.0$ s alongside a pure FEM and a pure MPM model. The hybrid FEM-MPM and pure MPM models yield very similar runout results, with a smooth final geometry sloping at an inclination approximate to the friction angle. This result is consistent with the simulation results of Lajeunesse (2005). On average, the final runout surface of the hybrid FEM-MPM model is 6.9% different from that of the pure MPM model. The pure FEM phase shows a distinctly different surface profile with more vertical movement, less runout, and a sharp corner.

This geometry results from compatible displacements of the elements due to the Lagrangian formulation and mesh distortion/entanglement at large deformations in FEM analysis.

These results confirm that MPM is more capable than FEM at capturing runout behavior, and the agreement between the results from the pure MPM and hybrid FEM-MPM analyses demonstrates the viability of the FEM to MPM transfer algorithm. The collapse of the granular column has a simple failure initiation mechanism – gravity-driven failure – which is why we do not see a significant difference between pure MPM and our hybrid FEM-MPM schemes. However, for problems where accurate quantification of stresses is essential for describing the failure kinematics, such as liquefaction-induced flow failures, the hybrid FEM-MPM approach will leverage the strengths of both FEM (failure initialization) and MPM (runout kinematics) methods.

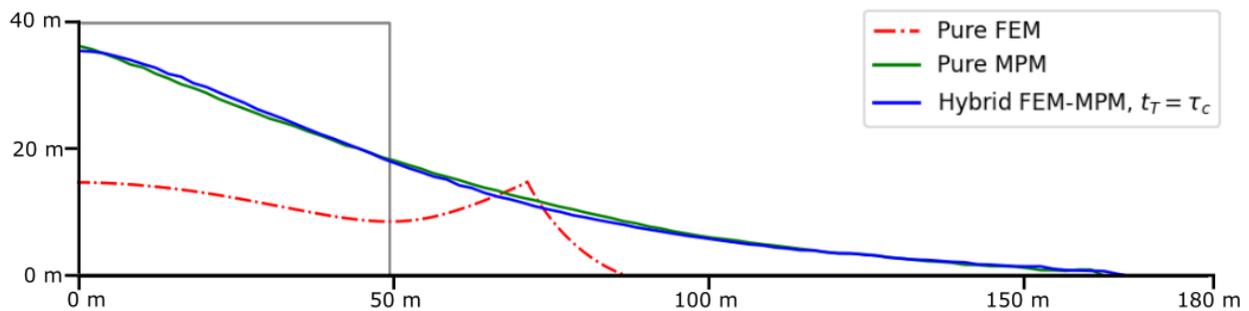

*Figure 6: Runout surfaces of pure FEM, pure MPM, and hybrid FEM-MPM models ($t_T$ = 2.0 s) of column collapse. While the pure MPM and hybrid results are final runout geometries, the pure FEM result is the surface at t = 8.0 s after the mesh has entangled.*

### *Effect of Transfer Time*

In the above example, we transferred the system from FEM to MPM at the critical time, which is well defined in a granular collapse experiment as the time at which the flow is fully mobilized. However, in most geotechnical slope failures, there is no theoretical critical time for failure mobilization, and there is no ideal transition time to transfer from FEM to MPM. In this section, we evaluate the effect of transfer time on the runout behavior of granular columns.

The runout behavior in granular columns is a function of its aspect ratio ($a$) and is independent of the geometric size of the problem. Short columns ($a < 2$) show a linear relation between runout and the initial aspect ratio, and their runout is characterized by frictional dissipation along a shear plane. In contrast, the runout distance for tall columns ($a > 2$) shows a power law relation with the initial aspect ratio, and these columns experience free-fall followed by collisional dissipation. We evaluate the effect of transfer times for a short column ($a = 0.8$), as well as a tall column ($a = 80$ m / $20$ m $= 4.0$). We consider these two columns because their failure mechanisms are different, so the effects of transfer time may vary.

Figure 7 shows the deviatoric strain ($\varepsilon_q$) distribution from the FEM phase of the models at the critical time of each model, as well as the time evolution of the surface profiles until mesh entanglement. For the short column (Figure 7a), we observe a diagonal failure plane above which the material slides. As time progresses, the surface profile for the short column displays the compatibility constraints in FE results in unrealistic geometries like upper right corner of the model, which influences the profile along the top and at the toe (Figure 7b). In contrast, for the tall column (Figure 7c), the entire volume above the failure surface experiences free fall. We do not observe unrealistic geometries, like the corner effect, for tall columns as their behavior is mainly controlled by the free fall mechanism (Figure 7d), and the geometry of the toe moves to the right, beginning to run out.

23/46

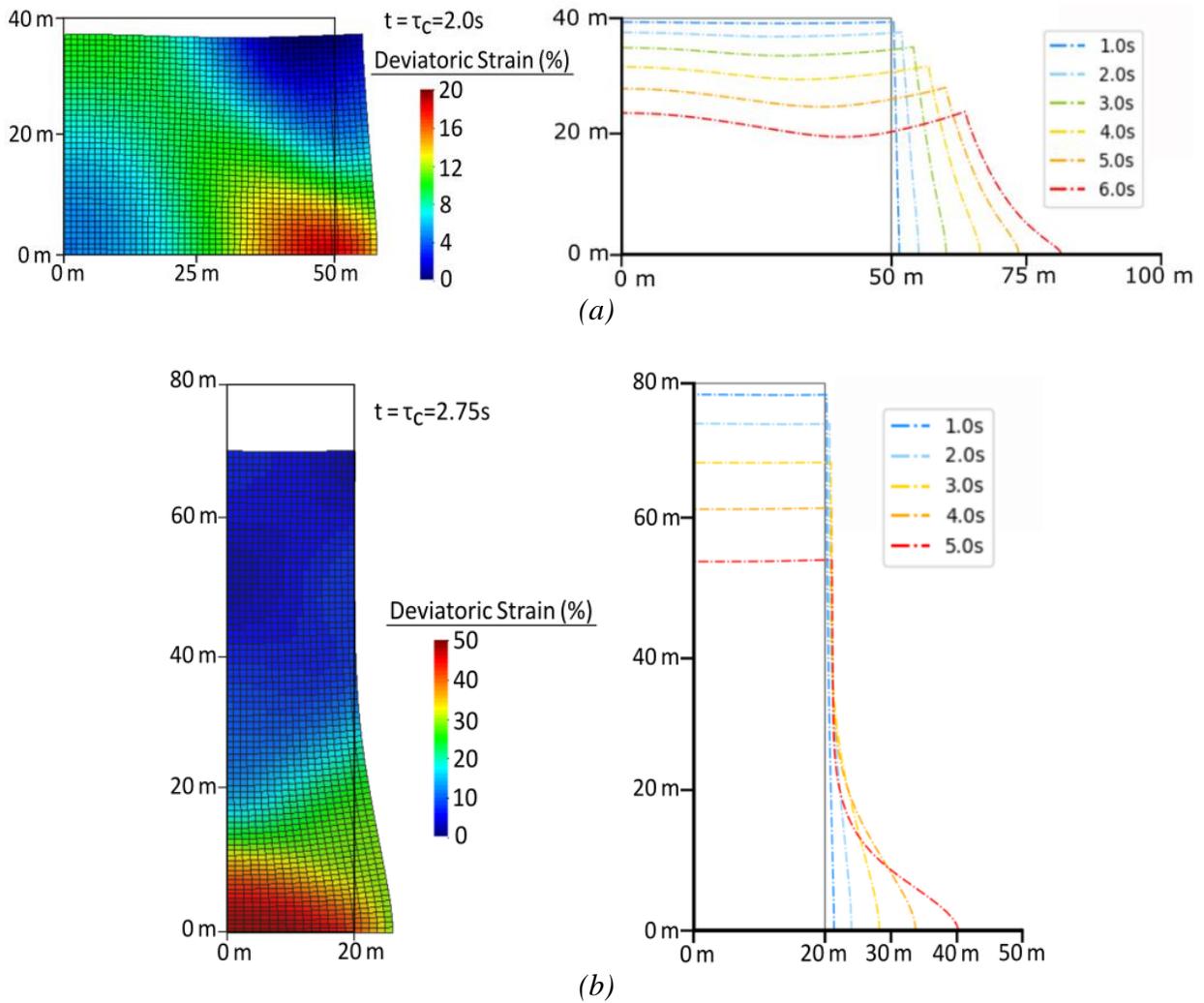

*Figure 7: FEM deviatoric strain ($\varepsilon_q$) contours of columns at $t = \tau_c$ and time evolution of FEM column shape until mesh entanglement for (a) short column (a = 0.8; $\tau_c$ ~ 2.0 s) and (b) tall column (a = 4.0; $\tau_c$ ~ 2.9 s).*

The key issue with the FEM as failure progresses is distortion and entanglement of the mesh at large deformations. Figure 8a shows the FEM mesh for the short column after entanglement at $t = 8$ s ($t/\tau_c$ ~ 2.75), when the top of the model has settled more than 20 m, and the runout ($L-L_o$) is 46 m. Since the boundary conditions are only applied to the nodes along the model's base, at large deformations many elements extend below the base of the model and fold over each other. Clearly, this level of mesh distortion is not realistic, and thus a transfer time of 8.0 s, or later, is not acceptable for this model.

To monitor the quality of the FEM mesh over time, we consider the determinant of the Jacobian for each element at time $t$, $|J_t|$, normalized to its initial value at the start of the simulation, $|J_0|$. As the runout evolves, the normalized Jacobian determinant, $|J_t|/|J_0|$, decreases, and these smaller values correspond to increased element distortion. A negative $|J_t|/|J_0|$ indicates that an element has become inverted and the mesh has entangled. Figure 8a shows the values of $|J_t|/|J_0|$ for each element at $t = 8.0$ s. Negative values of $|J_t|/|J_0|$ are observed in the elements that have folded over, while values less than 0.5 are observed within the diagonal shear zone. Rest of the mesh have $|J_t|/|J_0|$ greater than 0.5. Figure 8b tracks the minimum $|J_t|/|J_0|$ across all the FEM elements, representing the most distorted element. The rate of decrease for $(|J_t|/|J_0|)_{min}$ is small over the first ~2.0 s, but increases rapidly as the simulation progresses, suggesting that the quality of the FEM phase degrades quickly. Mesh entanglement, as indicated by $(|J_t|/|J_0|)_{min} < 0.0$, occurs in the short column just before $t = 6.0$ s ($t/\tau_c \sim 3.0$), cessation phase, and in the tall column just before $t = 5.0$ s ($t/\tau_c \sim 1.7$), during the runout phase.

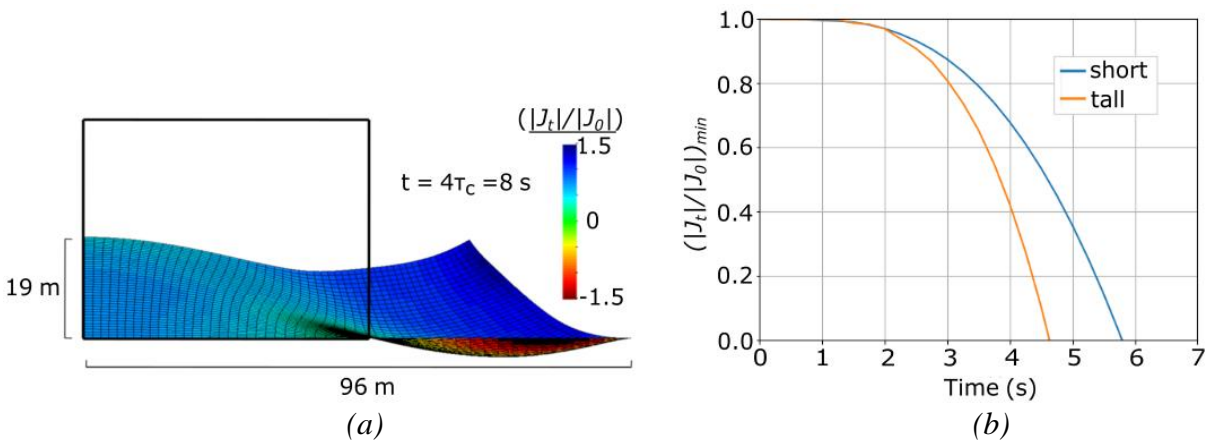

Figure 8: (a) Deformed mesh of the short column at $t = 8.0$ s with contours of $(|J_t|/|J_0|)$ and (b) $(|J_t|/|J_0|)_{min}$ vs. time for both columns.

We select a range of transfer times ($t_T$) for the FEM phase of column collapse to evaluate the effect of $t_T$ on the runout kinematics of the hybrid FEM-MPM analysis. To avoid the effects of mesh entanglement, we only consider transfer times that are associated with the $(|J_t|/|J_0|)_{min} > 0.0$ (i.e., $t_T < 6.0$ s for the short column, $t_T < 5.0$ s for the tall column). We perform transfers at 1.0 s intervals. A transfer at 0.0 s means that the transfer occurs immediately after the linear elastic geo-static phase and uses a plastic constitutive model after the transfer. Each hybrid FEM-MPM model continues until the column settles into its final configuration.

The surface profile results of the hybrid FEM-MPM models are shown in Figure 9, along with the results from a pure MPM model. For both columns, the earliest transfers resemble the pure MPM models most closely, while later transfers become increasingly different. Transfers shortly after 0.0 s yield slightly larger runout distances, while later transfer times yield smaller runouts. For the short column, the transfer time also noticeably influences the final crest height, with later transfer times yielding more crest settlement. Significant crest settlement occurs during the FEM phase (Figure 7b), so later transfers inherit this behavior. Regardless of transfer time, all hybrid FEM-MPM models show a smooth runout surface like the pure MPM model, unlike the FEM results that generate a distorted and unrealistic runout profile (Figure 7). However, the variation in final runout distances and crest settlement among the hybrid FEM-MPM models indicates that the transfer time from FEM to MPM impacts the runout results.

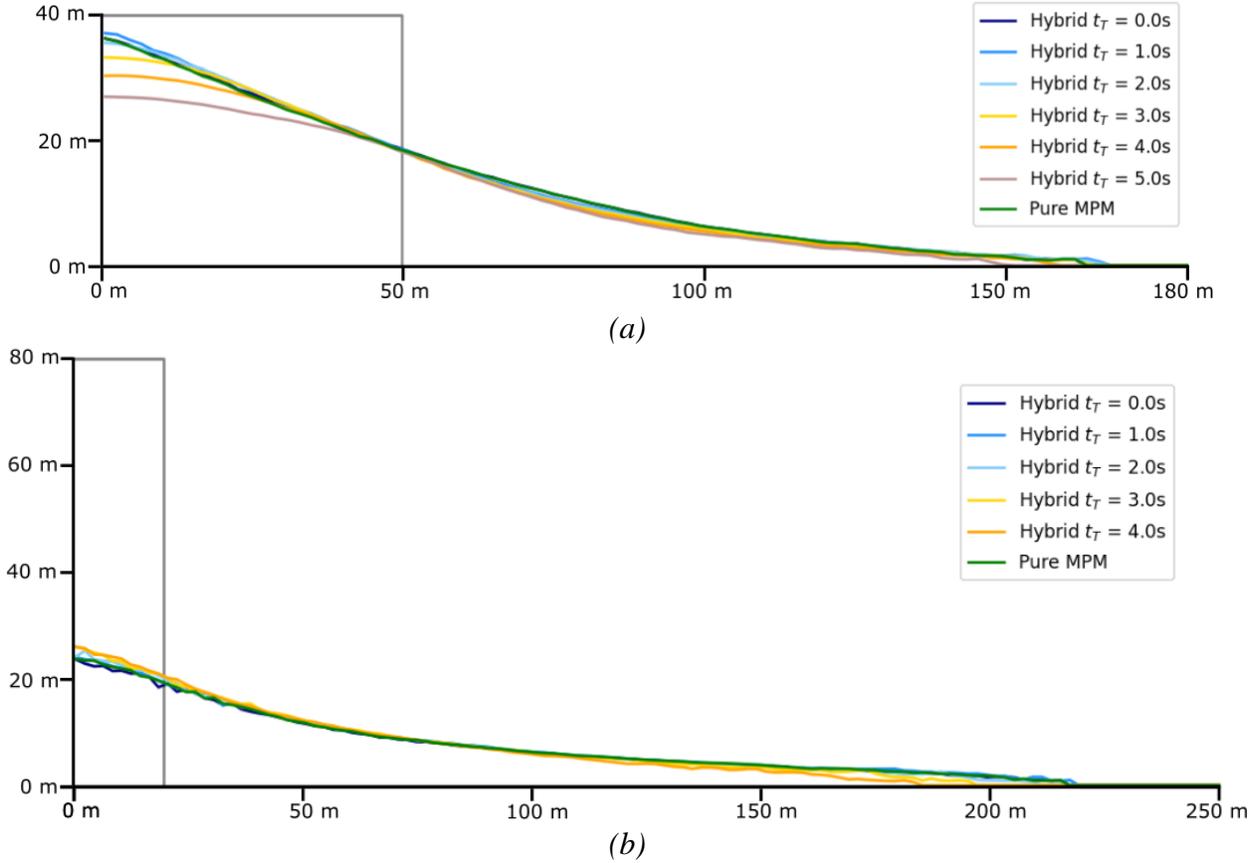

*Figure 9: Final runout surface profiles from hybrid FEM-MPM models and pure MPM model for (a) short column and (b) tall column collapse models.*

To quantify the error in the final geometries of the different models, we define three error metrics: the relative runout error ($R_N$ Error), the relative crest height error ($H_c$ Error), and the root mean square error (RMSE) of the full runout profile relative to the $t_T = 0.0$ s hybrid FEM-MPM model. We compute errors relative to the $t_T = 0.0$ s hybrid FEM-MPM model because our goal is to evaluate errors associated with different transfer times. The runout profiles are defined, approximately, by the 99th percentile y-coordinate of all the material points within a 2-m wide horizontal bin. We select 99th percentile to avoid extraneous particles. The RMSE is defined from the runout profiles using:

$$RMSE = \sqrt{\frac{\sum_{b=1}^{N_{bins}}\left(\frac{y_b - y_{b,0}}{y_{b,0}+p}\right)^2}{N_{bins}}} \qquad (17)$$

where $N_{bins}$ is the total number of bins sampled, $y_b$ is the elevation of bin $b$ in the model, $y_{b,0}$ is the elevation of bin $b$ in the sequential hybrid model for a transfer time $t_T = 0.0$ s, and $p$ is the particle size. The particle size is included to avoid excessively large errors near the toe where the elevations are small. For the RMSE calculation, we consider bins from the left boundary to the toe of the $t_T = 0.0$ s model.

The surface profile error metrics are compiled in Table 2, and the errors increase as the transfer time gets larger. The RMSE metric generally quantifies the largest error because this metric is influenced significantly by the differences in height near the toe where the values of $y_{b,0}$ are small. For the latest transfer times, which are close to the time of mesh entanglement, the RMSE is greater than 30% and the runout and crest height errors are greater than 10%. The larger errors occur because the later transfers inherit the unrealistic displacements and mesh distortions from the FEM phase, which affect the MPM runout results.

Also shown in Table 2 are the surface profile errors of the pure MPM models. These errors are all less than ~ 3%, indicating that the small errors associated with the RBF stress transfer do not significantly influence the runout results.

*Table 2: Relative errors of final geometrics of hybrid FEM-MPM column models with different transfer times.*

| Column | $t_T$ (s) | $t_T/\tau_c$ | $R$ (m) | $R_N$ | $R_N$ Error | $H_c$ | $H_c$ Error | RMSE |
|---|---|---|---|---|---|---|---|---|
| 50x40 | 0 | 0.00 | 110.5 | 2.21 | -- | 36.2 | -- | -- |
| | 1 | 0.50 | 113.8 | 2.28 | 3.0% | 37.1 | 2.5% | 8.1% |
| | 2 | 0.99 | 114.4 | 2.29 | 3.5% | 35.5 | -1.9% | 8.1% |
| | 3 | 1.49 | 111.2 | 2.22 | 0.6% | 33.2 | -8.3% | 7.1% |
| | 4 | 1.98 | 106.2 | 2.12 | -3.9% | 30.3 | -16.3% | 18.7% |
| | 5 | 2.48 | 97.8 | 1.96 | -11.5% | 26.9 | -25.7% | 31.2% |
| | Pure MPM | | 110.7 | 2.21 | 0.2% | 36.2 | 0.0% | 3.0% |
| 20x80 | 0 | 0.00 | 194.8 | 9.74 | -- | 23.3 | -- | -- |
| | 1 | 0.35 | 197.7 | 9.89 | 1.5% | 24.0 | 3.0% | 10.3% |
| | 2 | 0.70 | 188.0 | 9.40 | -3.5% | 25.2 | 8.2% | 16.0% |
| | 3 | 1.05 | 175.5 | 8.78 | -9.9% | 26.0 | 11.6% | 29.2% |
| | 4 | 1.40 | 163.4 | 8.17 | -16.1% | 26.0 | 11.6% | 39.5% |
| | Pure MPM | | 195.6 | 9.78 | 0.4% | 23.7 | 1.7% | 4.4% |

To investigate the progression of the runout of the models, we consider the time evolution of the normalized runout ($R_N$) and of the total kinetic energy ($KE$) normalized by the potential energy of the initial state of the FEM phase ($PE_0$). For the MPM phase, we calculate the kinetic energy as:

$$KE^{MPM} = \sum_{i=1}^{N_p} \frac{1}{2} m_i^p v_i^{p^2} \qquad (18)$$

where $N_p$ is the number of particles in the model, $m_i^p$ is the mass of the particle, and $v_i^p$ is the velocity of the particle. For the FEM phase, we calculate the kinetic energy of the model as:

$$KE^{FEM} = \sum_{i=1}^{N_e} \frac{1}{2} m_i^e v_i^{e^2} \qquad (19)$$

where $N_e$ is the number of elements in the model, $m_i^e$ is the mass of the element, and $v_i^e$ is the average velocity of the four nodes connected to the element. The potential energy of the initial state of the FEM phase ($PE_0$) is calculated as:

$$PE_0 = \sum_{i=1}^{N_e} m_i^e g h_i^e \qquad (20)$$

where $N_e$ is the number of elements in the model, $g$ is the acceleration of gravity, and $h_i^e$ is the average elevation of the four nodes connected to the element.

Figure 10 displays the normalized runout ($R_N$) and normalized kinetic energy ($KE/PE_0$) as a function of normalized time ($t/\tau_c$) for the different column models and transfer times. The runout progressions of all the hybrid FEM-MPM models generally display the theoretical three-phase pattern of the granular column collapse, as previously illustrated in Figure 4, and the peak in $KE/PE_0$ occurs at the transition from the initiation phase to runout phase. The early transfers match the theoretical time evolution more closely, with the runout phase initiating at $t/\tau_c \sim 1$, which also corresponds with the time of the peak $KE/PE_0$. For the later transfer times, the runout phase is delayed because the Lagrangian nature of the FEM phase slows down the runout initiation. For these later transfer times, there is an acceleration in $KE/PE_0$ after the transition but the peak $KE/PE_0$ is smaller and the peak occurs later. The smaller $KE/PE_0$ is due to the later transfers inherit unrealistic displacements and distortion from the FEM phase, as shown by the FEM profiles in Figure 7. This effect is illustrated most significantly for the $t_T = 5.0$ s short column model, which experiences a sudden decrease in $KE/PE_0$ after the transfer (Figure 10a) and represents the most distorted short column and the lowest $(|J_t|/|J_0|)_{min}$ at the time of transfer (Figures 7 and 8b). As a result of the smaller $KE/PE_0$ for the later transfer times, these models experience smaller runout.

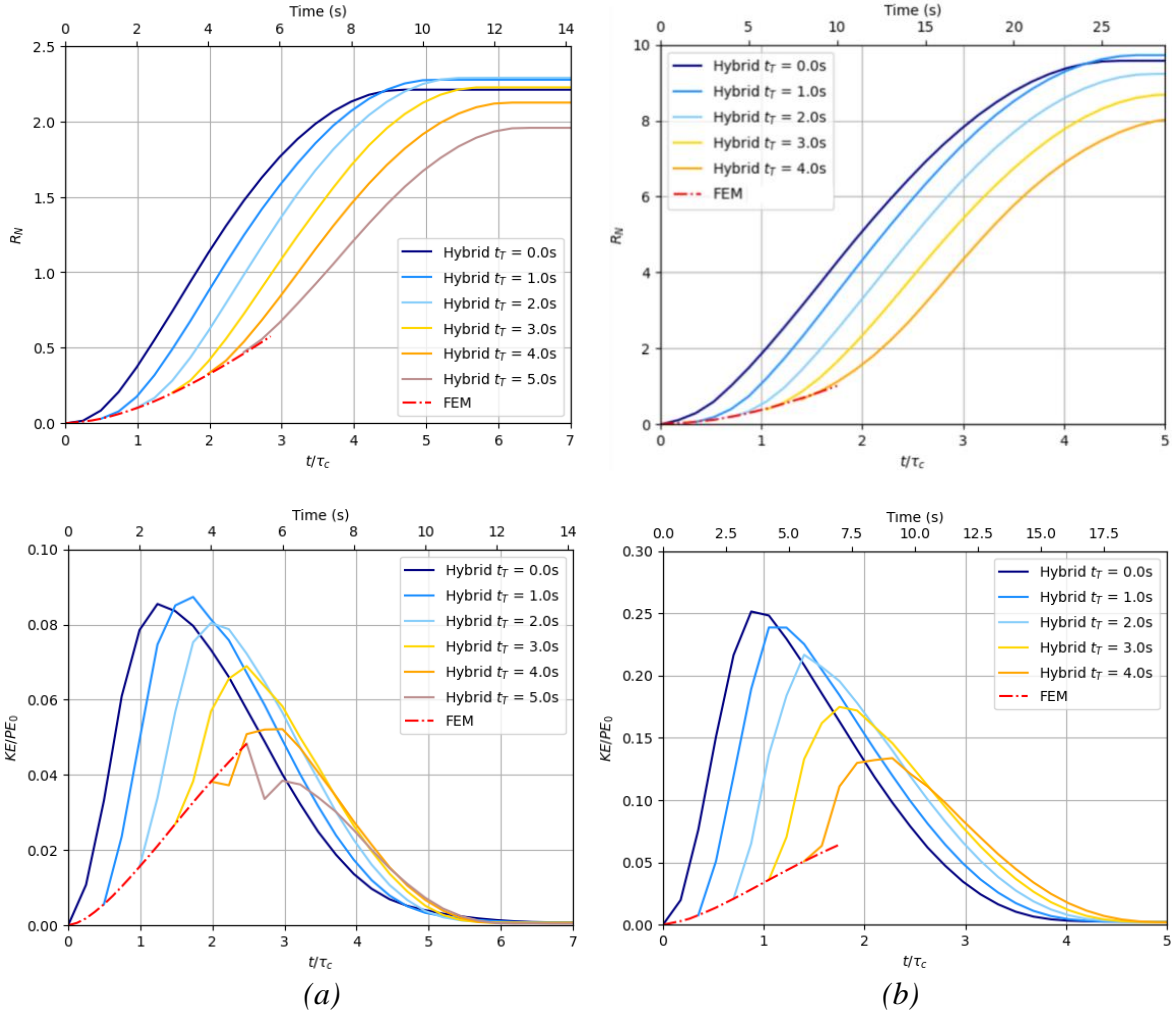

*Figure 10: Time evolution of normalized runout ($R_N$) and normalized kinetic energy ($KE/PE_0$) for (a) short column and (b) tall column models.*

These hybrid FEM-MPM analyses of column collapse illustrate that the runout response of the MPM phase is significantly influenced by the displacement and mesh distortion present in the FEM phase at the time of transfer. This issue will be important when considering complex failure mechanisms in which the FEM phase must be utilized as long as possible to capture failure initiation. In such cases, the transfer will need to be performed at an intermediate time after the FEM phase has developed the failure initiation but before its mesh has become excessively distorted.

# EVALUATION OF SEQUENTIAL HYBRID FEM-MPM APPROACH FOR SLOPE FAILURES

## *Slope Failures*

We next evaluate the hybrid FEM-MPM approach for a more common geotechnical problem: earth slope failures and subsequent runout. For this evaluation, we consider two slope geometries: a 20 m high slope inclined at 1H:1V consisting of friction-dominated soil and a 10 m high slope inclined at 3H:1V consisting of cohesion-dominated soil (Figure 11). We evaluate the performance of the hybrid FEM-MPM method on the two different slope models because they will have distinct failure characteristics: a shallow failure surface for the frictional slope and a deep failure surface for the cohesive slope.

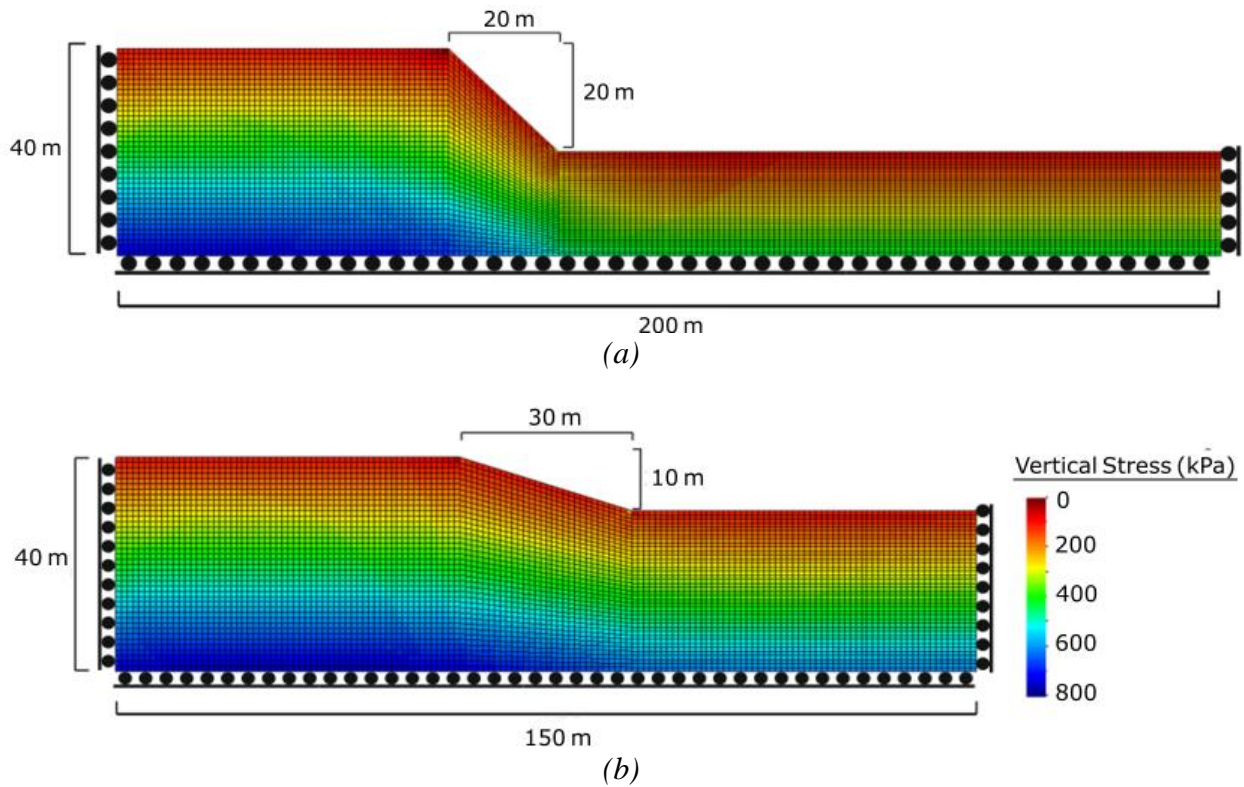

*Figure 11: FEM elastic vertical stress for a) frictional slope model and b) cohesive slope model (stress in kPa and negative indicates compression).*

Similar to the column collapse experiments, we initialize the FEM analysis of the slope models by computing the geo-static stresses with a linear elastic constitutive model, and failure is initiated by gravity by changing to a plastic constitutive model. The strength parameters (Table 3) are assigned to achieve a factor of safety of 0.5 in both slopes, as computed using the stability charts of Wright (1969). We analyze each slope geometry with the hybrid FEM-MPM approach, considering different transfer times to assess the effect of transfer time on the kinematic behavior and runout predictions.

*Table 3: Properties of frictional and cohesive slope material.*

| Linear Elastic (both slopes) | | |
|---|---:|---|
| Young's Modulus | 172 | Mpa |
| Poisson's Ratio | 0.23 | |
| Frictional Slope | | |
| Friction Angle | 10 | degrees |
| Cohesion | 0.5 | kPa |
| Mass Density | 1925 | kg/m3 |
| Cohesive Slope | | |
| Friction Angle | 3.6 | degrees |
| Cohesion | 5.8 | kPa |
| Mass Density | 1690 | kg/m3 |

Figure 12 shows the $\varepsilon_q$ contours of the FEM phases at $t$ = 1.0 s for the frictional slope and at $t$ = 2.0 s for the cohesive slope, times at which the failures have begun to mobilize. Also shown are the FEM surface profiles at different times. The $\varepsilon_q$ values within the failure surfaces range from about 8 to 16%, with the largest strains localized at the toe for the frictional slope (Figure 12a) but extending over a larger section of the failure surface for the cohesive slope (Figure 12b). As a result, the FEM phase of the cohesive slope reaches larger displacements at similar magnitudes of strain. For the meshes shown, $(|J_t|\,/\,|J_0|)_{min}$ is equal to 0.98 for both slopes. The FEM meshes entangle (i.e., $(|J_t|\,/\,|J_0|)_{min} < 0.0$) at $t \sim 5.0$ s for the frictional slope

and $t \sim 11.0$ s for the cohesive slope. The surface profiles of the frictional slope are mostly smooth (Figure 12a), except at the toe, until $t \sim 5.0$ s when significant toe bulging erupts. The surface profiles of the cohesive slope are also initially smooth (Figure 13b), but toe bulging becomes significant by $t \sim 5.0$ s and by $t \sim 10.0$ the toe of the slope extends higher than the original slope height. These FEM surface profiles at later times are unrealistic and represent further evidence that an FEM mesh is incapable of accurately capturing the large deformations and runout associated with slope failure.

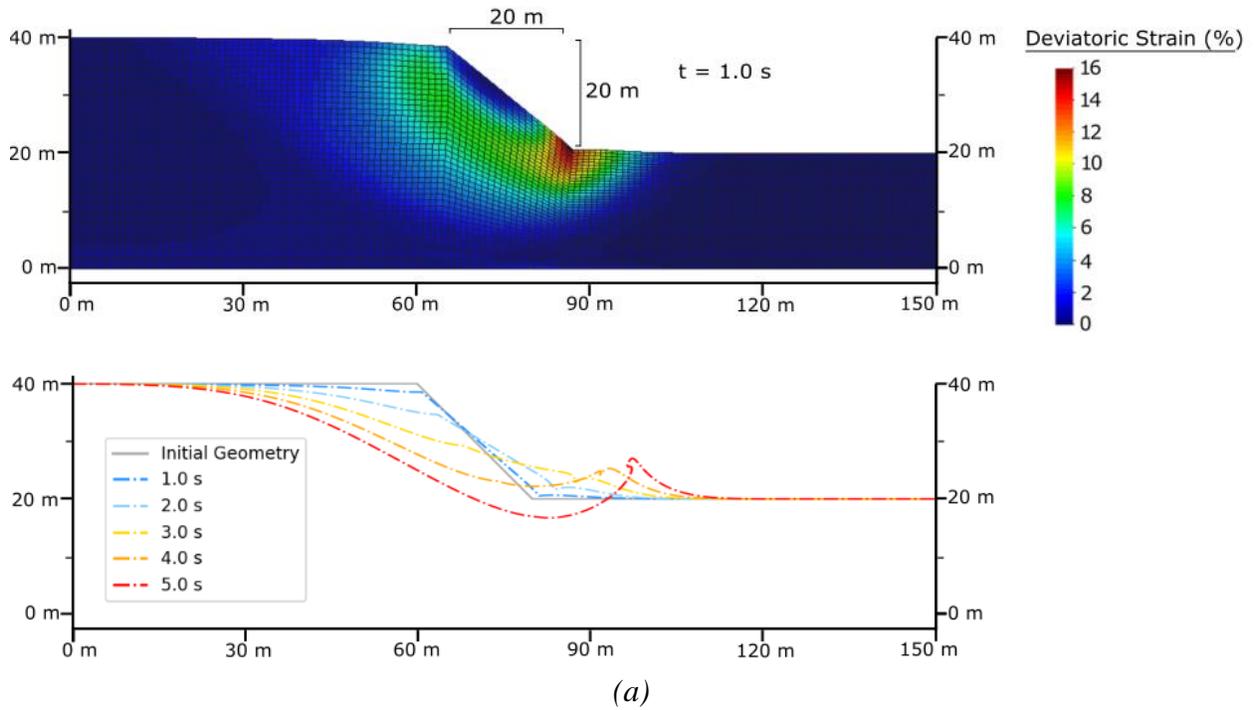

(a)

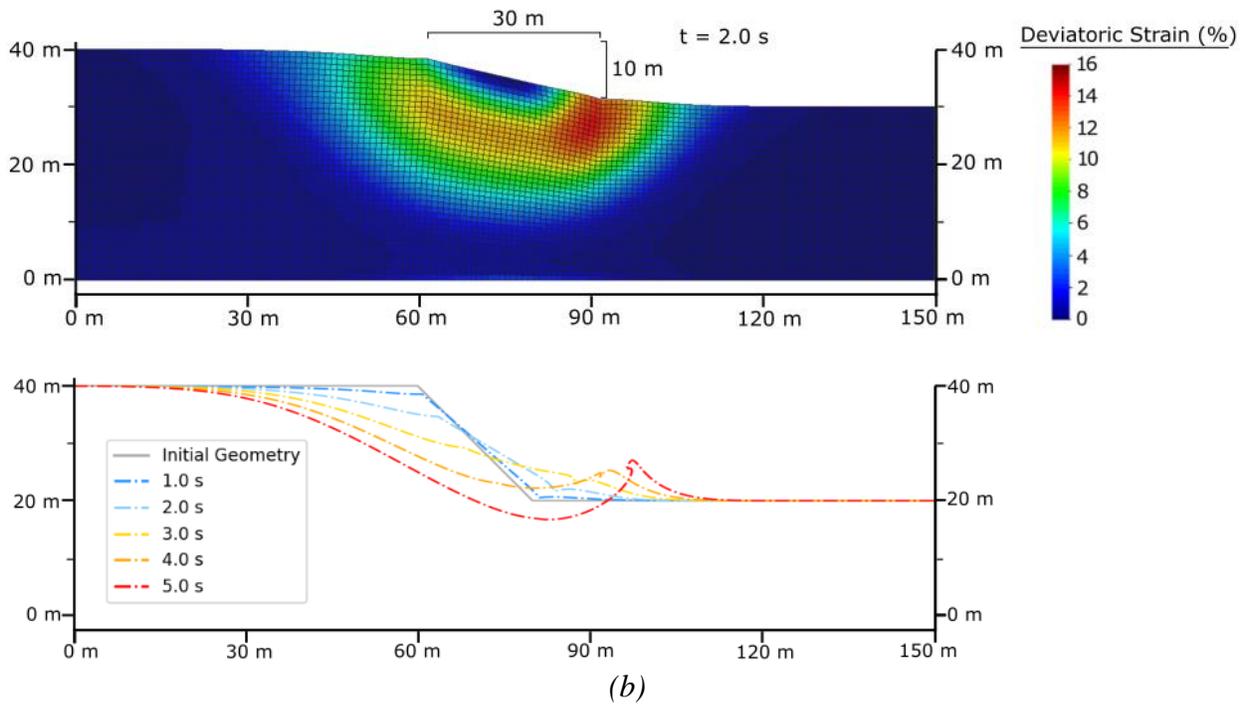

Figure 12: Deviatoric strain ($\varepsilon_q$) contours and surface profiles for the (a) frictional slope and (b) cohesive slope.

## *Effect of Transfer Time*

For each slope model, a set of $t_T$ values is selected between $t = 0$ s and the time of mesh entanglement, and the FEM phase is transferred to MPM at each $t_T$ and allowed to run until movement stops. The FEM phases all have one Gauss point per element and each element is transferred into sixteen material points. The MPM phase of these analyses sometimes experience material points detaching from the model at the flow front, so we apply 1% Cundall damping to to minimize this behavior. However, even with the additional damping, the $t_T = 10.0$ s transfer of the cohesive slope was numerically unstable.

Figure 13 shows the final surface profiles of the hybrid FEM-MPM models for each $t_T$. For the frictional slope, the runout decreases with increasing $t_T$ (Figure 13a), similar to the analyses of granular column collapse. This similarity is expected given that both the fictional slopes and granular columns feature frictional materials with similar failure mechanisms. For the cohesive slope, there is clearly an effect of $t_T$ on the final surface profile (Figure 13b), but the deeper failure surface does not develop a distinct runout location. The later transfer times ($t_T$ = 5.0 and 7.5 s) display more crest settlement and more toe bulging than the $t_T$ = 0.0 s model, because the severe crest settlement and toe bulging from the FEM phase (Figure 12b) is inherited by the MPM phase.

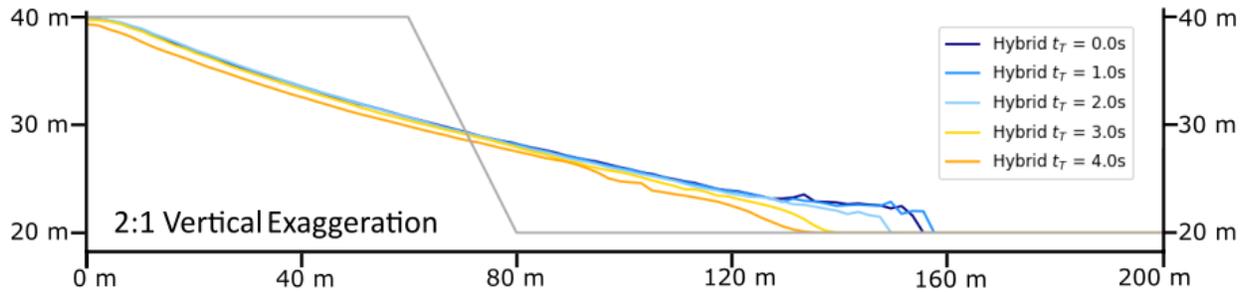

*(a)*

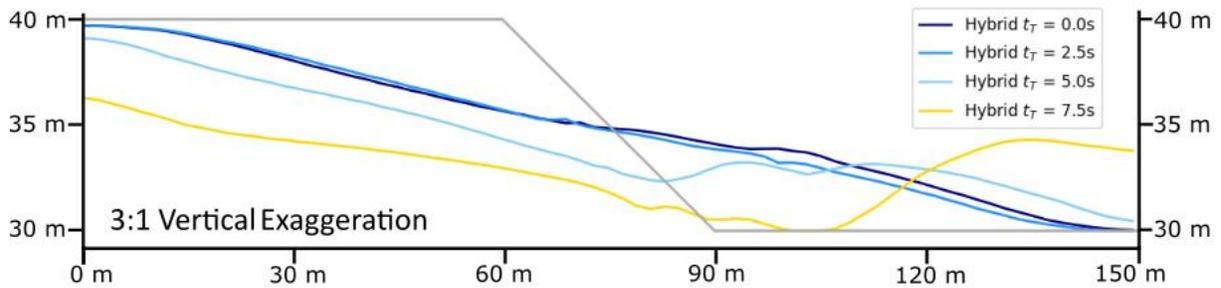

*(b)*

*Figure 13: Final surface profiles for all hybrid FEM-MPM models of (a) frictional and (b) cohesive slope models.*

Table 4 tabulates the runout error and the RMSE of the surface profiles for the slopes as a function of transfer time. For these simple, gravity-driven failures, we again compute the error relative to the $t_T$ = 0.0 s model. Because the cohesive slopes do not have a distinct runout location, the runout error is not defined for these results. For the frictional slope, the runout error and RMSE increase with increasing $t_T$, with values as large as 20% to 40% for $t_T \geq 3$ s. The RMSE for the cohesive slope also increases with increasing $t_T$, but the values of RMSE are much larger than for the frictional slope and reach almost 400% for $t_T = 7.5$ s. These results indicate that the cohesive slope suffers more than the frictional slope from the later transfer times, revealing that the significance of the transfer time may be different for different failure mechanisms.

*Table 4: Runout and RMSE error of hybrid slope models.*

| Slope | $t_T$ (s) | Runout (m) | R Error | RMSE |
|---|---|---|---|---|
| Frictional | 0 | 74 | -- | -- |
| | 1 | 74 | 0.0% | 5.6% |
| | 2 | 69 | -6.8% | 19.3% |
| | 3 | 59 | -20.3% | 35.1% |
| | 4 | 54 | -27.1% | 41.5% |
| Cohesive | 0 | -- | -- | -- |
| | 2.5 | -- | -- | 14.1% |
| | 5 | -- | -- | 68.1% |
| | 7.5 | -- | -- | 384.0% |
| | 10 | -- | -- | unstable |

To evaluate the kinematics of the hybrid FEM-MPM models of the slope failures, Figure 14 plots the time evolution of $KE/PE_0$ for different $t_T$. For the frictional slope, the progressions of $KE/PE_0$ for the earlier transfers (i.e., $t_T \leq 2.0$ s) all display an increase in $KE/PE_0$, relative to the FEM phase, after the transfer. The $KE/PE_0$ progression for the FEM phase shows a similar shape to the MPM phase of these hybrid models, but it reaches its peak at a later time. For the MPM hybrid models, the peak $KE/PE_0$ is smaller for the later transfer times, which results in smaller runouts. The $t_T = 3.0$ s and 4.0 s hybrid models display a different response with a sudden deceleration immediately after the transfer, which results in a more severe reduction in runout than for the earlier transfers (Figure 13; Table 4). The meshes at these later transfer times are more distorted with $(|J_t|/|J_0|)_{min}$ of 0.53 and 0.03 and $\varepsilon_{q,max}$ of 85% and 123%, respectively.

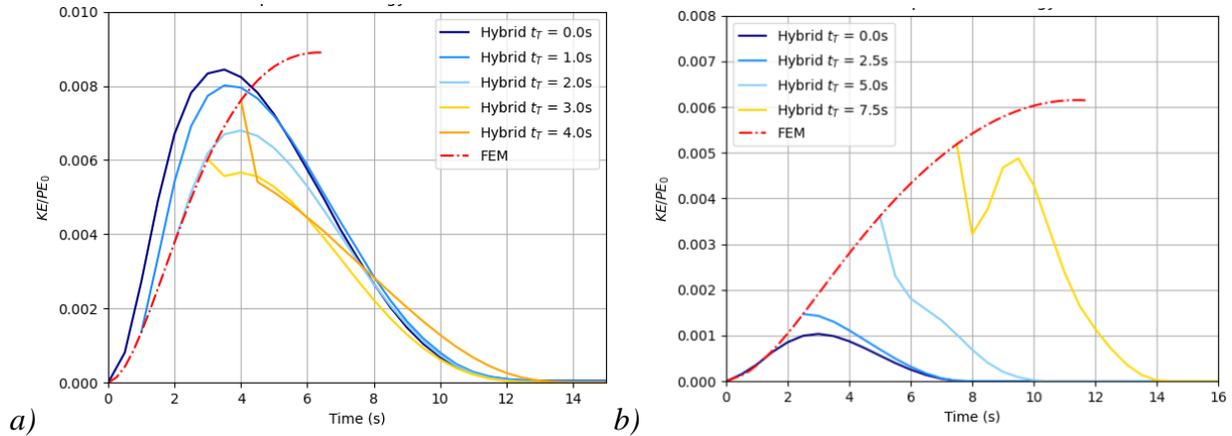

*Figure 14: Progression of normalized kinetic energy in a) frictional slope model and b) cohesive slope model.*

For the cohesive slope (Figure 14b), the $KE/PE_0$ for the $t_T = 0.0$ s model behaves similarly to the FEM phase until $t \sim 2$ s, but its $KE/PE_0$ then decelerates, reaching its peak at $t \sim 3$ s and reaching zero at $t \sim 7$ s. The $KE/PE_0$ of the FEM phase continues to increase with time and reaches its peak much later ($t \sim 11$ s) and with a much greater value of $KE/PE_0$. None of the transfers of the cohesive slopes experience the sudden increase in $KE/PE_0$ observed in the frictional slope. Rather, the $t_T = 2.5$ s model slowly decelerates, the $t_T = 5.0$ s model exhibits a sudden deceleration after the transfer, and the $t_T = 7.5$ s model uniquely experiences a period of deceleration and then a subsequent period of acceleration before finally settling towards zero. These responses are a direct result of mesh distortion in the FEM phase and its resulting unrealistic deformations prior transfer. By $t = 5.0$ s and 7.5 s, $(|J_t|/|J_0|)_{min}$ has fallen to 0.77 and 0.42, respectively, and $\varepsilon_{q,max}$ has reached 56% and 91%. The toe bulging in the FEM phase at $t = 7.5$ s is as large as 9 m, nearly the height of the original slope, which severely influences the final surface profile of the hybrid model results.

The results shown here further indicate that the selection of transfer time exerts a significant impact on the final surface profiles of the slope models. For the frictional slope, later transfers exhibit reduced runout distances and gradually increasing error, the same trend as observed for the column collapses. For the cohesive slope, the later transfers exhibit much greater RMSE than those of the frictional slope and the final surface profiles do not even qualitatively resemble the profile for the $t_T = 0.0$ s transfer. We attribute these differences to the distinct failure mechanism of the cohesive slope, where the failure surface for the cohesive slope is deeper and larger than that of the frictional slope, so it experiences larger displacements before the mesh entangles.

**RELATIONSHIP BETWEEN MESH DISTORTION AND RUNOUT ERROR**

In practice, choosing several values of $t_T$ and conducting MPM simulations for each may be impractical. A user will therefore need to decide when to transfer the model from FEM to MPM based only on the conditions of the FEM phase. In complex problems like those involving earthquake-induced liquefaction, $t_T = 0$ s will not be an option, as FEM is necessary to capture pore pressure generation during earthquake shaking and the development of the failure mechanism. Therefore, it will be necessary to allow some degree of mesh deformation before the transfer to MPM. While we have qualitatively shown that more mesh distortion leads to larger error in the runout, here we offer a quantitative relationship between the two.

We seek to relate the RMSE (Equation 17) of the slope surface of each hybrid FEM-MPM model to the quality of the FEM mesh at the time of the transfer. We quantify the quality of the mesh considering two metrics: $|J_t|/|J_0|$ and $\varepsilon_q$. These values are averaged across all elements with $\varepsilon_q > 3\%$, creating two values which indicate the degree of mesh distortion at a given time. Unsurprisingly, these two metrics of mesh distortion are highly correlated.

Figure 15 presents the RMSE for the results of the hybrid FEM-MPM models, along with the corresponding $\varepsilon_{q,avg}$ and $(|J_t|/|J_0|)_{avg}$ of the FEM phase at $t_T$. The data presented includes the short (50 m x 40 m) and tall (20 m x 80 m) columns as well as the frictional slope presented in this work. These data are supplemented with data from multiple transfers of a 20 m x 20 m column, an 80 m x 80 m column, and a frictional slope with a factor of safety is 0.7. The cohesive slope model is excluded as it deviates heavily from the trends observed in the models with frictional soils (e.g., the failure mode is fundamentally different; the RMSEs are much higher). Therefore, the data in Figure 15 are only representative of failures in frictional materials. The data indicate that RMSE increases with increasing mesh distortion, as quantified by increasing $\varepsilon_{q,avg}$ and decreasing $(|J_t|/|J_0|)_{avg}$. The slope models tend to yield slightly larger RMSE than the column models with the same degree of mesh distortion, but the trend is relatively consistent across the different models.

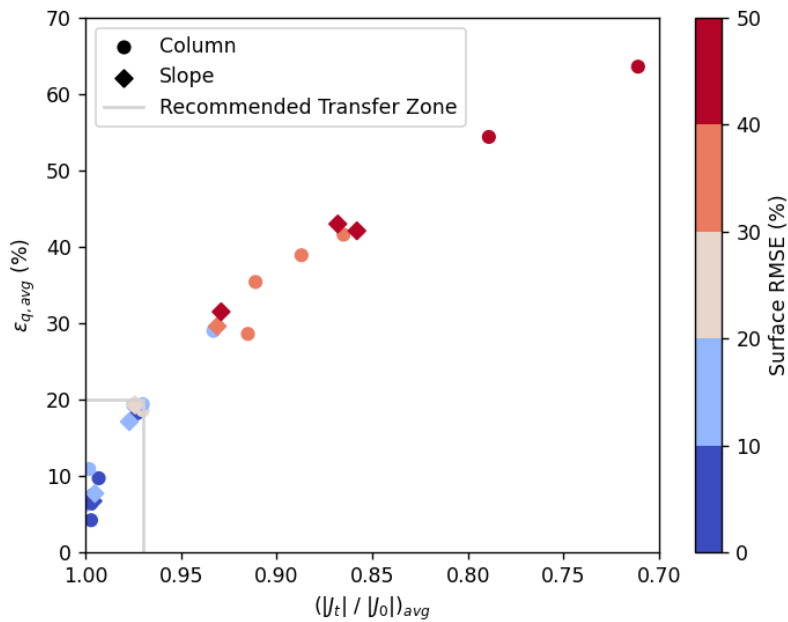

*Figure 15: Surface profile error (RMSE) as a function of mesh quality ($\varepsilon_{q,avg}$ and $(|J_t| / |J_0|)_{avg}$ across elements with $\varepsilon_q > 3\%$) of frictional models.*

We consider an RMSE of 30% to be a threshold below which the results are reasonable and usable for engineering evaluation. The latest transfer times which yield hybrid models that meet this accuracy criterion are $t_T = 4.0$ s for the short column, $t_T = 3.0$ s for the tall column, and $t_T = 2.0$ s for the frictional slope. Transfers performed before $(|J_t| / |J_0|)_{avg}$ has fallen below 0.97 and before $\varepsilon_{q,avg}$ reaches 20% consistently deliver RMSE values below 30%. This range of mesh quality is shown in Figure 15 as the "Recommended Transfer Zone". Therefore, we recommend that, as far as possible, transfers be performed before the mesh distorts beyond these levels to ensure that the results of the model remain reasonable and useful.

**CONCLUSIONS**

The FEM is a well-established numerical modeling technique that excels in predicting the conditions of landslide initiation but is incapable of modeling large displacements and runouts, while MPM is capable of accounting for large deformations but lacks the stress precision and boundary conditions of FEM for predicting failure initiation. We propose a novel hybrid FEM-MPM approach for slope failure modeling that sequentially employs the FEM for failure initiation and then switches to MPM for the runout. Through this hybrid approach, we utilize the strengths of each technique to predict both landslide initiation and runout with a single model of maximal accuracy. We evaluate the proposed hybrid approach by simulating granular column collapse experiments and simple slope failures to identify a range of suitable times to transfer from FEM to MPM.

The transfer from FEM to MPM is performed at a user-specified time. The final runout results can be significantly impacted by the transfer time, so it is critical to understand the restraints to a suitable transfer time. The accuracy of the FEM phase decreases with mesh distortion, and any inaccuracy accumulated in the FEM phase is inherited by the MPM phase and eventually the final results. Therefore, the earliest possible transfer time should be selected to minimize mesh distortion. In the simple models evaluated here, the transfer can occur at $t_T = 0.0$ s, but in more complex cases, especially those involving seismic triggers and liquefaction, a later transfer time will be required because only the FEM phase is capable of accurately capturing the triggering effects of the earthquake (i.e. the distribution of liquefaction). We found that, in failures of frictional soils, transfers performed before $(|J_t| / |J_0|)_{avg}$ and $\varepsilon_{q,avg}$ among elements with $\varepsilon_q > 3\%$ fall below/exceed 0.95 and 20%, respectively, avoid severe error. The criteria for suitable transfer in failures involving cohesive soils are less clear, but our experiments indicate that severe error is preceded by qualitatively unrealistic geometries in the FEM phase. The ideal transfer time will thus be case-dependent and range between the full development of the failure mechanism in FEM and the development of excessive mesh deformation.